\documentclass[a4paper]{iopart}
\usepackage{iopams}
\usepackage{graphicx}
\usepackage{psfrag}
\usepackage{dcolumn}
\usepackage{slashbox,multirow}
\usepackage{amsthm,amssymb}
\usepackage{cases,color}
\usepackage{mathrsfs}
\usepackage{times}
\DeclareSymbolFontAlphabet{\mathrsfs}{rsfs} \DeclareMathAlphabet{\mathcal}{OMS}{cmsy}{m}{n}
\newcommand{\scri}{\mathrsfs{I}}

\newcommand{\be}{\begin{equation}}
\newcommand{\ee}{\end{equation}}
\renewcommand{\figurename}{Figure}

\begin{document}

\title{Universality of global dynamics for the cubic wave equation}

\author{Piotr Bizo\'n$^{1}$ and An{\i}l Zengino\u{g}lu$^{2,3}$}

\address{$^{1}$ M. Smoluchowski Institute of
 Physics, Jagiellonian University, Krak\'ow, Poland}

\address{$^{2}$ Max-Planck-Institut f\"ur Gravitationsphysik
 (Albert-Einstein-Institut), Golm, Germany}

\address{$^{3}$ Department of Physics, and Center for Scientific
 Computation and Mathematical Modeling, University of Maryland,
 College Park, MD 20742, USA}

\date{\today}
\begin{abstract}
We consider the initial value problem for the spherically symmetric, focusing cubic wave equation
in three spatial dimensions. We give numerical and analytical evidence for the existence of a
universal attractor which encompasses both global and blowup solutions. As a byproduct we get an
explicit description of the critical behavior at the threshold of blowup.
\end{abstract}


\section{Introduction} We consider a  semilinear wave equation in three spatial dimensions with
the focusing cubic nonlinearity
\begin{equation}\label{eqo}
\partial_{tt} \phi - \Delta\phi- \phi^3=0\,.
\end{equation}
Heuristically, the dynamics of solutions of this equation can be viewed as a competition between the
Laplacian  which tends to disperse the waves and the nonlinearity which tends to concentrate the
waves. For small initial data the dispersion "wins" leading to global solutions which decay to
zero as $t\rightarrow \infty$. For large initial data the dispersive spreading is too weak to
counterbalance the focusing nonlinearity and solutions blow up in  finite time. For each of these
generic evolutions the leading asymptotic behavior  is known: small solutions decay as $1/t^2$ at
timelike infinity \cite{ch1,nik}, while large solutions diverge as $\sqrt{2}/(T-t)$ for $t$
approaching a blowup time $T$ \cite{mz1,bct}. The dichotomy of dispersion and blowup brings up
the question of what determines the boundary between these two behaviors and, in particular, what
is the evolution of critical initial data which lie on the boundary.
During numerical investigations of this question we observed  that for a
large set of initial data the solutions rapidly converge to a universal attractor which is given
by a two-parameter family of explicit solutions of equation (\ref{eqo}). The aim of this paper is
to give analytic and numerical evidence for this behavior (which is rather surprising for a
conservative wave equation). The description of the critical dynamics, which motivated our
investigations, emerges as a special case.

The paper is organized as follows. In section~2 we recall some basic properties of solutions of
equation (\ref{eqo}) and we formulate three conjectures about the existence of the attractor.
Analytic evidence for these conjectures is given in section~3. After a short description of our
numerical methods and tests in section~4, we present numerical evidence for the conjectures in
section~5. Finally, in section~6 we make some general remarks.

\section{Preliminaries and conjectures}
In this paper we restrict our attention to spherically symmetric solutions $\phi=\phi(t,r)$, so
equation (\ref{eqo}) reduces to
\begin{equation}\label{eq}
\partial_{tt} \phi - \partial_{rr}\phi-\frac{2}{r}\partial_r \phi - \phi^3=0\,.
\end{equation}
This equation  is invariant under the following transformations:
\begin{itemize}
\item translation in time $T_a$ by a constant $a$
\begin{equation}\label{trans}
T_a:    \phi(t,r)\, \rightarrow \, \phi(t+a,r)\,,
\end{equation}
\item scaling $S_{\lambda}$ by a positive constant $\lambda$
\begin{equation}\label{scal}
S_\lambda:    \phi(t,r)\, \rightarrow \, \frac{1}{\lambda}
   \phi\left(\frac{t}{\lambda},\frac{r}{\lambda}\right)\,,
\end{equation}
\item conformal inversion $I$ (which is an involution)
\begin{equation}\label{inv}
I:    \phi(t,r)\, \rightarrow \, \frac{1}{t^2-r^2} \, \phi\left(\frac{t}{r^2-t^2},
   \frac{r}{t^2-r^2}\right)\,,
\end{equation}
\item reflection $\phi \rightarrow -\phi$.
\end{itemize}
\noindent Neglecting the Laplacian in (\ref{eqo}) and solving the ordinary differential equation
$\partial_{tt}\phi-\phi^3=0$, one gets the one-parameter family of spatially homogeneous solutions
\begin{equation}\label{u0}
   \phi_b=\frac{\sqrt{2}}{b-t}\,.
\end{equation}
This family is a special case of the two-parameter family of solutions \cite{anco}
\begin{equation}\label{orbit}
   \phi_{(a,b)}(t,r)=\frac{\sqrt{2}}{t+a+b\left((t+a)^2-r^2\right)}\,,
\end{equation}
which can be obtained from (\ref{u0}) by the action of conformal inversion $I$ followed by the
time translation $T_a$.
The scaling transformation $S_{\lambda}$ acting on
$\phi_{(a,b)}$ only rescales the parameters of the solution without changing its form.
Note
that the solution $\phi_{(a,b)}(t,r)$ is singular on the two-sheeted hyperboloid $t= -a -1/(2b)\pm
\sqrt{1/(4b^2)+r^2}$.

The main result of this paper is the observation that the family (\ref{orbit}) is an attractor
for a large set of initial data. More precisely, we have the following conjectures about the
forward-in-time behavior of solutions of equation (\ref{eq}) starting from smooth, compactly
supported (or suitably localized) initial data (by time reflection symmetry analogous conjectures
can be formulated for the backward-in-time behavior):
\vskip 0.2cm \noindent{\emph{Conjecture~1.}} For any generic globally regular solution
$\phi(t,r)$ there exist parameters $(a,b)\in \mathbb{R}\times \mathbb{R}^+$ and $\kappa=\pm 1$
such that $t^4 \left(\phi(t,r)-\kappa\, \phi_{(a,b)}(t,r)\right)$ is bounded for all finite $r$
and $t\rightarrow \infty$.
\vskip 0.2cm \noindent{\emph{Conjecture~2.}}
For any solution $\phi(t,r)$ which blows up at the origin in finite time $T$,
 there exist parameters $(a,b)\in \mathbb{R}\times \mathbb{R}^-$ and $\kappa=\pm 1$ such that
   $(T-t)^{-2} \left(\phi(t,r)-\kappa\, \phi_{(a,b)}(t,r)\right)$
is bounded for $t\rightarrow T^-$ inside the past light cone of the blowup point.
\vskip 0.2cm \noindent{\emph{Conjecture~3.}} The borderline between dispersive and blowup
solutions, described respectively in Conjectures~1 and~2, consists of codimension-one globally
regular solutions $\phi(t,r)$ for which there exist a parameter $a\in \mathbb{R}$ and $\kappa=\pm
1$ such that
   $t^4 \left(\phi(t,r)-\kappa\, \phi_{(a,0)}(t,r)\right)$
is bounded for all finite $r$ and $t\rightarrow \infty$. \vskip 0.2cm \noindent The parameters
$(a,b)$ for which the above assertions hold will be referred to as optimal. \vskip 0.2cm

\noindent \emph{Remark~1:} Conjectures 1 and 2 are refinements of the well-known asymptotic
behavior of solutions of equation (\ref{eq}): $t^{-2}$ decay near timelike infinity in the case
of global regularity (for small initial data) and $\sqrt{2}/(T-t)$ growth in the case of blowup
(for large initial data), respectively.

\noindent \emph{Remark 2:} Although Conjecture~3 is a special limiting case
of Conjecture~1, we state it as a separate conjecture to emphasize the critical character
of the attractor solution $\phi_{(a,0)}$. Note that the slowly decaying global solutions
described in Conjecture~3 are not asymptotically free, i.e., they do not scatter.

\noindent \emph{Remark~3:} The genericity condition in Conjecture~1 is essential
because, as we shall see below, there do exist non-generic very rapidly
decaying globally regular solutions which do not converge to the attractor $\phi_{(a,b)}$.

\noindent \emph{Remark~4:} For very large amplitudes the blowup first occurs on a sphere
\cite{bct} and then Conjecture~2 does not hold.

In the remainder of the paper we give evidence for the above conjectures. The evidence is based
mainly on numerical simulations, however, before presenting numerics, we will give two analytic
arguments: one based on linearized stability analysis and one based on an explicit solution.

\section{Analytic evidence}\label{sec:3}
\subsection{Linearized stability}
In this section we discuss the linearized stability of the attractor solution
$\phi_{(a,b)}(t,r)$. Since linearization commutes with symmetries, we may set $a=b=0$ without
loss of generality. We denote the resulting solution by $\phi_0$, thus $\phi_0=\sqrt{2}/t$. In
order to determine the spectrum of small perturbations around this forward self-similar solution
we will use the symmetry under conformal inversions. To this end, consider the region $I^+(O)$,
the interior of the future light cone of the origin $(0,0)$ of the Minkowski spacetime, and its
foliation by spacelike hyperboloids (in this paragraph we
 follow Christodoulou \cite{ch2})
\begin{equation}\label{hyp}
   t=c+\sqrt{c^2+r^2}\,,
\end{equation}
where $c$ is a positive constant. The hyperboloid (\ref{hyp}) is asymptotic to $\partial I^+(c)$,
the future light cone of the point $(c,0)$. Let us make the conformal inversion
\begin{equation}\label{ci}
   I: (t,r) \mapsto (\bar t,\bar
   r)=\left(\frac{t}{r^2-t^2},\frac{r}{t^2-r^2}\right)\,.
\end{equation}
This transformation maps $I^+(O)$ (in the original coordinate system) to $I^-(O)$ (the interior
of the past light cone of the origin in the barred coordinate system). In particular, the
hyperboloids (\ref{hyp}) are mapped to spacelike hyperplanes
\begin{equation}\label{hplanes}
   \bar t = - \bar c\,, \qquad \mbox{where} \qquad \bar c=\frac{1}{2c}\,.
\end{equation}

Now, consider a global-in-time solution $\phi(t,r)$ inside $I^+(O)$. In the barred coordinates we
have from (\ref{inv})
\begin{equation}\label{phib}
   \bar \phi(\bar t,\bar r)=(t^2-r^2) \phi(t,r)\,.
\end{equation}
It follows from the above paragraph that the study of the asymptotics of $\phi(t,r)$ for
$t\rightarrow \infty$ is equivalent to the study of the asymptotics of $\bar\phi(\bar t,\bar r)$
for $\bar t \rightarrow 0^-$. In particular, there is equivalence between linearized
perturbations of $\phi_0=\sqrt{2}/t$ for $t\rightarrow \infty$ and linearized perturbations of
$\bar \phi_0=\sqrt{2}/(-\bar t)$ for $\bar t\rightarrow 0^-$. But the latter have already been
determined in \cite{bct,gp}. Namely, it has been shown there that the spectrum of smooth linearized
perturbations about the backward self-similar solution $\sqrt{2}/(-\bar t)$ consists of a
discrete set of eigenmodes of the form $(-\bar t)^{\lambda_n} \xi_n(y)$, where
$y=\bar r/(-\bar t)$ and $\lambda_0=-2,\lambda_1=0,\lambda_n=n \,(n\geq 2)$. More precisely, we have the
following eigenmode expansion around $\bar \phi_0$
\begin{equation}\label{bct}
 \!\!\!\!\!\!\!\!\!\!\!\!\!\!\!\!\!\!\!\!\!\!\!\!\!\!  \delta\bar \phi(\bar t,\bar r)=c_0 (-\bar t)^{-2}
   +c_1 (1-y^2) +c_2 (-\bar t)^{2}
   (1-\frac{2}{3}y^2+\frac{1}{5}y^4)+
   \mathcal{O}((-\bar t)^3)\,.
\end{equation}
By (\ref{phib}), this translates into the following eigenmode expansion around $\phi_0$ (where
now $y=r/t$)
\begin{equation}\label{exp2}
\!\!\!\!\!\! \delta\phi(t,r)= c_0 (1-y^2)+c_1 \frac{1}{t^2} +
   c_2 \frac{1}{t^4} \frac{1-2y^2/3+y^4/5}{(1-y^2)^3}+
   \mathcal{O}(1/t^5)\,.
\end{equation}
The first two eigenmodes in the expansions (\ref{bct}) and (\ref{exp2}) correspond to the perturbations
along the symmetry
orbits $\bar \phi_{(a,b)}(\bar t,\bar r)$ and $\phi_{(a,b)}(t,r)$ respectively. For example, we have
\begin{equation}\label{gena}
\!\!\!\!\!\!\!\!\! \frac{\partial}{\partial b}
   \phi_{(a,b)}(t,r) \vert_{(a=0,b=0)}
   \sim 1-y^2\,,\quad
   \frac{\partial}{\partial a} \phi_{(a,b)}(t,r)
   \vert_{(a=0,b=0)} \sim \frac{1}{t^2}\,.
\end{equation}
The choice of the optimal parameters $(a,b)$ for the attractor amounts to tuning away the
coefficients of the symmetry modes $c_0$ and $c_1$, hence the rate of convergence to the
attractor is expected to be governed by the third eigenmode in the expansions (\ref{bct}) and
(\ref{exp2}):
\begin{equation}\label{conver0}
   \bar \phi(\bar t,\bar r)-\bar \phi_{(a,b)}(\bar t,\bar r) \sim (-\bar t)^2 \qquad \mbox{for}
   \,\,\bar t\rightarrow 0^-\,,
\end{equation}
and
\begin{equation}\label{conver}
   \phi(t,r)-\phi_{(a,b)}(t,r) \sim 1/t^4 \qquad \mbox{for}
   \,\,t\rightarrow \infty\,,
\end{equation}
which is consistent with our conjectures. Below we will verify this expectation numerically, but
first we want to give a simple example which corroborates  (\ref{conver}).
\subsection{The conformal solution}
Equation (\ref{eq}) has the following globally regular explicit solution
\begin{equation}\label{A2}
   \phi_{conf}(t,r)= \frac{2}{\sqrt{(1+(t-r)^2)(1+(t+r)^2)}}\,,
\end{equation}
which we will refer to  as the conformal solution. Note that this solution corresponds to
time symmetric initial data
\begin{equation}\label{ic}
   \phi(0,r)=\frac{2}{1+r^2}\,,\qquad \partial_t \phi(0,r)=0\,.
\end{equation}
The conformal solution can be easily found in the framework of conformal compactification which
maps Minkowski space with the flat metric $\eta$ into the Einstein universe with the conformal
metric $g=\Omega^{2} \eta$, where the conformal factor $\Omega=2[(1+(t-r)^2)(1+(t+r)^2)]^{-1/2}$
\cite{Penrose65}.
Under this mapping the cubic wave equation $\Box \phi+\phi^3=0$ transforms to $\Box_{g} \Phi
-\Phi+\Phi^3=0$, where $\Phi=\Omega^{-1}\phi$. The trivial constant solution $\Phi=1$ gives
$\phi_{conf}=\Omega$.

By elementary calculation we find that for the conformal solution  the optimal parameters of the
attractor $\phi_{(a,b)}$ are $(a,b)=(-1/\sqrt{2},1/\sqrt{2})$. More precisely, we have for
$t\rightarrow\infty$
\begin{equation}\label{a2v}
   \phi_{conf}(t,r)-\phi_{(-\frac{1}{\sqrt{2}},\frac{1}{\sqrt{2}})}(t,r)
   = -\frac{1}{t^4}\, \frac{3+y^2}{(1-y^2)^3} +
   \mathcal{O}(t^{-6})\,,
\end{equation}
in agreement with Conjecture~1.


\section{The numerical method and tests}
\subsection{The hyperboloidal initial value problem for the cubic wave equation}
\begin{figure}[t]
 \begin{centering}
   {\psfrag{t}{$t$} \psfrag{r}{$r$}
     \includegraphics[width=0.31\textwidth,height=0.17\textheight]{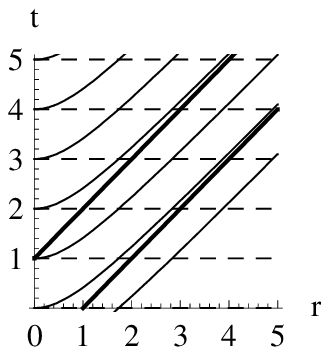}
     \hspace{23mm}} {\psfrag{ip}{$i^+$} \psfrag{i0}{$i^0$}
     \psfrag{scrp}{$\scri^+$}
     \includegraphics[width=0.27\textwidth]{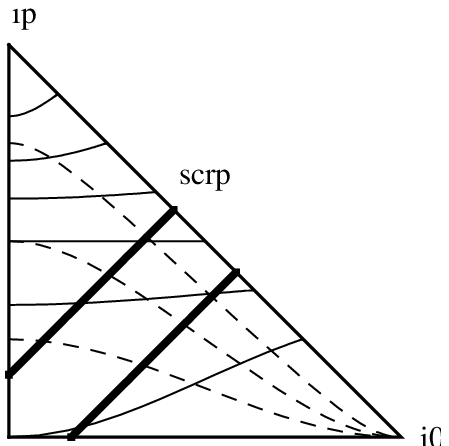}} \caption{
     The future domain of a $t=0$ surface is partially depicted on
     the left panel in standard coordinates and entirely represented
     on the right panel in a Penrose diagram
     \cite{Penrose65,Hawking73,Zeng07}.  Dashed lines are level sets of $t$,
     solid lines are hyperboloids shifted in time as given by
     (\ref{eq:tau}), the thick straight lines depict outgoing
     characteristics indicating the location of the approximate
     wavefront. \label{fig:1}}
 \end{centering}
\end{figure}

In order to study the asymptotic behavior of solutions and verify the conjectured convergence to
the attractor, it is convenient to foliate Minkowski spacetime by hyperboloids. A time coordinate
$\tau$ adapted to such a foliation can be written as \be\label{eq:tau} \tau = t -
\sqrt{\frac{9}{K^2}+r^2}.\ee The level surfaces of $\tau$ are standard hyperboloids shifted in
the time direction. They have constant mean curvature $K$ which is a free parameter. In the
following we set $K=3$ for simplicity. A foliation of Minkowski spacetime by level sets of $\tau$
is depicted in Figure \ref{fig:1}.

In order to be able to analyze the propagation of the outgoing waves to infinity we
introduce a compactifying radial coordinate $\rho$ along the surfaces of our foliation. For the
regularity of our equations in this compactified setting, we perform a suitable conformal
rescaling of the metric. The rescaling factor denoted by $\Omega$ must be a function of $\rho^2$
to ensure regularity at the origin in the conformal manifold. We choose $\Omega=(1-\rho^2)/2$
following \cite{Moncrief00, Husa02, Fodor03}. The compactifying coordinate is then chosen
according to $r=\rho/\Omega$.

We express the standard Minkowski metric $\eta$ in the new coordinates $(\tau,\rho)$ and rescale
it with the conformal factor $\Omega^2$ to obtain \be\label{eq:metric} g = \Omega^2\eta =
-\Omega^2\,d\tau^2 - 2\rho\,d\tau d\rho + d\rho^2 + \rho^2\,d\sigma^2, \ee where $d\sigma^2$ is
the standard metric on the unit two-sphere. In these coordinates the zero set of the conformal factor corresponds to future null infinity denoted by $\scri^+$. We rewrite the cubic wave equation on Minkowski
spacetime in conformally covariant form \be \label{eq:conf} \Box_g \Phi -\frac{1}{6} R[g]\, \Phi
+ \Phi^3 = 0, \quad \mathrm{where} \quad \Phi = \frac{\phi}{\Omega}\,, \quad g=\Omega^2 \eta\,.
\ee
The Ricci scalar of the metric $g$ from (\ref{eq:metric}) appearing in the above equation is given by
\[R[g]=\frac{12\, (1-\rho^2)\, (3+\rho^2)}{(1+\rho^2)^3}\,.\]
With the auxiliary variables,
\[
\psi:=
\partial_\rho \Phi\quad \mbox{and}
\quad \pi := \frac{2}{1+\rho^2}(\partial_\tau\Phi+\rho\,\partial_\rho \Phi)\,,\]
we can rewrite the
system (\ref{eq:conf}) in first order symmetric hyperbolic form as
\begin{eqnarray}\label{eq:wave}
\partial_\tau \Phi &=&\frac{1+\rho^2}{2}\,\pi - \rho\,\psi, \nonumber
\\ \partial_\tau \psi &=&\partial_\rho\left(\frac{1+\rho^2}{2}\pi -
\rho\,\psi\right), \\ \partial_\tau \pi &=&
\frac{1}{\rho^2}\,\partial_\rho\left(\rho^2\left(\frac{1+\rho^2}{2}
\psi-\rho\,\pi\right)\right) +
\frac{1+\rho^2}{2}\left(\Phi^3-\frac{1}{6}\,R[g] \Phi \right)
\nonumber.
\end{eqnarray}

The initial data are arbitrary in our conjectures. In the numerical studies presented below we
choose a Gaussian pulse \be\label{eq:id} \Phi(0,\rho) = A\,e^{-(\rho-\rho_c)^2/\sigma^2}, \qquad
\partial_\tau \Phi(0,\rho) = 0, \ee  with fixed  parameters $\rho_c=0.3,~\sigma=0.07$, and
varying amplitude $A$. Tests for different initial data give the same qualitative results which
makes us feel confident that the phenomena described below are universal.
\subsection{The code}
We solve the hyperboloidal initial value problem (\ref{eq:wave},\ref{eq:id}) numerically
using 4th order Runge-Kutta integration in time
and 6th order finite differencing in space. At the origin we apply the regularity condition
$\psi(\tau,0) = 0$. No boundary conditions are needed at the outer boundary because there are no
incoming characteristics due to the compactification. One-sided finite differencing is applied on
the numerical boundaries.

\begin{figure}[t]
 \begin{centering}
   \psfrag{t}{$\tau$} \psfrag{Q}{$Q$}
   \includegraphics[width=0.42\textwidth]{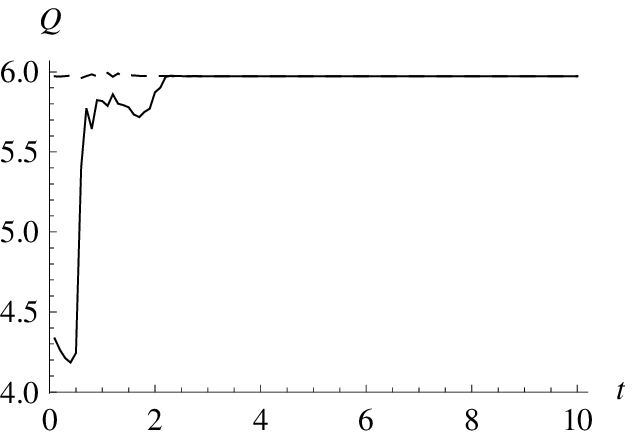}\hspace{15mm}
   \includegraphics[width=0.42\textwidth]{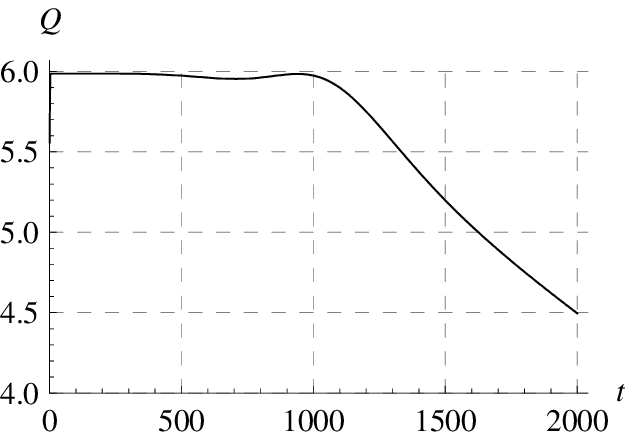}
   \caption{Convergence factors in time measured in the
     $L_2$-norm. The convergence factor $Q$ is defined by
     \mbox{$Q=\log_2\frac{ \| \Phi^{low}-\Phi^{med}\|}{\|
         \Phi^{med}-\Phi^{high}\|}$}. On the left panel we see that
     after a short transient phase the code converges with 6th
     order. The dashed curve corresponds to a simulation with a 10
     times smaller Courant factor and shows 6th order convergence
     from the beginning. The right panel shows results of a long time
     convergence test and indicates loss of convergence after about
     $\tau=1000$. \label{fig:conv}}
 \end{centering}
\end{figure}

To test the code we performed a three level convergence study with
200, 400 and 800 grid cells on the coordinate domain $\rho\in [0,1]$
that has infinite physical extent. For this study, we used an
amplitude of $A=2$ which is below the critical amplitude, and a
Courant factor of $\triangle \tau/\triangle \rho=0.8$. The convergence
factors shown in Figure \ref{fig:conv} indicate that, as expected, our
code is 6th order convergent after a short transient phase. The dashed
curve has been calculated using a much smaller Courant factor of
$0.08$ to show that the initial transient phase is due to numerical
errors in the time integration that converge with 4th order. The long
time convergence plot on the right panel in Figure \ref{fig:conv}
shows that convergence is lost at about $\tau=1000$ for the number of
cells used in this study. This is due to accumulation of numerical
errors and occurs at later times in tests with higher resolution. We
give numerical evidence for universal dynamics only in the convergent
regime. Two-level convergence tests using the explicit conformal
solution (\ref{A2}) give the same qualitative results.

\begin{figure}[t]
 \begin{centering}
   {\psfrag{t}{$\tau$} \psfrag{p}{$p_\rho$}
     \includegraphics[width=0.41\textwidth]{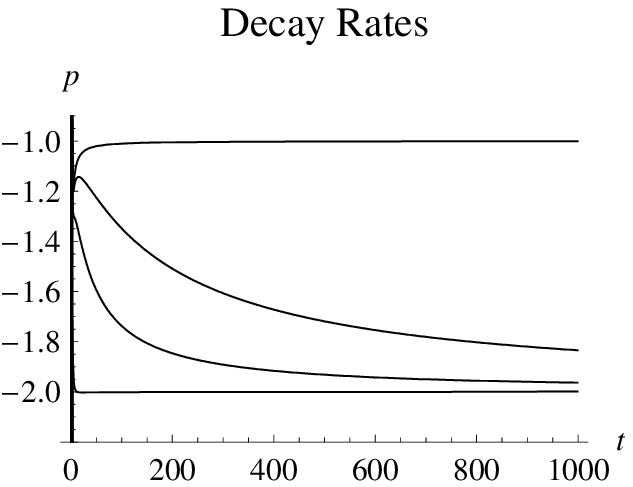}\hspace{13mm}}
   {\psfrag{t}{$\frac{1}{(T-\tau)}$}\psfrag{phi}{$\Phi$}
     \includegraphics[width=0.41\textwidth]{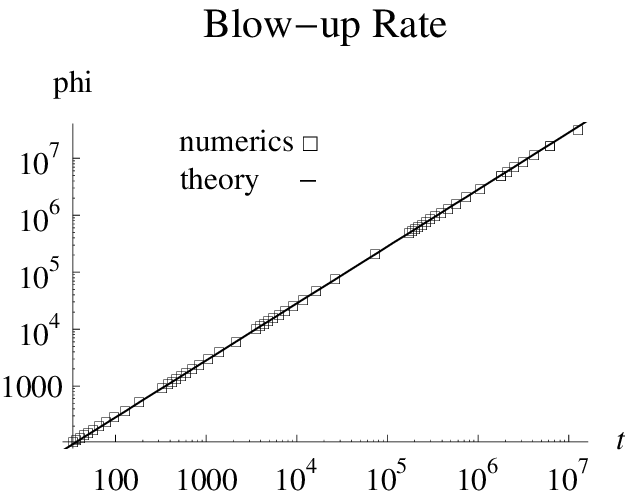}}
   \caption{Left panel: The local power index calculated along
     $\scri^+$, \mbox{$r=\{100,20\}$} and at the origin from top
     to bottom. Right panel: A numerical solution that blows up at
     the origin is depicted by small squares and the theoretical
     prediction of the blowup is depicted by the solid line on a
     log-log scale. The small squares are not distributed uniformly
     in time because we reduce our time steps while the solution
     grows. The figures indicate that the code reproduces the known
     decay rates for small solutions and the blowup rate for large
     solutions very accurately. \label{fig:tests}}
 \end{centering}
\end{figure}

Further tests can be performed using known properties of solutions to the cubic wave equation as
studied in \cite{bct}. First, it is known that generic small solutions decay as $t^{-2}$ near
timelike infinity and $t^{-1}$ along null infinity \cite{nik}. It is clear from (\ref{orbit})
that the attractor solution has this behavior for $b>0$. In order to accurately measure the
decay rate of the solution along null infinity, near timelike infinity and in the transition
domain between these two asymptotic regimes, we calculate the local power index, $p_\rho(\tau)$,
defined by
\[ p_\rho(\tau) = \frac{d\ln |\Phi(\tau,\rho)|}{d\ln \tau}. \]
Figure \ref{fig:tests} shows that our code reproduces the predicted
decay rates very accurately.

Second, it is known that large initial data lead to blowup in finite
time with the blowup rate at the origin $\sqrt{2}/(T-t)$
\cite{mz1,bct}. This is in accordance with the attractor solution
(\ref{orbit}) for $b<0$. The right panel in Figure \ref{fig:tests}
shows the numerical solution at the origin for large data on a log-log
scale against the theoretical prediction of the blowup rate. The two
curves match over many orders of magnitude indicating that our code
can reliably handle the blowup.
\section{Numerical evidence for universal dynamics}
The space of solutions to the cubic wave equation can be divided into three parts: decay, blowup
and criticality. These parts correspond respectively to $b>0$, $b<0$, and $b=0$ for the attractor
solution (\ref{orbit}). We will follow this natural classification in our presentation of the
numerical evidence for the universality of dynamics.

In many cases, the numerical evidence will be presented in terms of the conformally rescaled
solution in the coordinates presented in the previous section. In these variables the
two-parameter family of solutions (\ref{orbit}) takes the form (using the abbreviation
$\tilde{\tau}=\tau+a$)
\begin{equation} \label{phiab_hypal} \!\!\!\!\!\!\!\Phi_{(a,b)}(\tau,\rho) =
\frac{2 \sqrt{2}}{(\tilde{\tau}+1) (b\, (\tilde{\tau}+1)+1)-\rho^2
 (\tilde{\tau}-1) (b\, (\tilde{\tau}-1)+1)}\,. \end{equation}

\subsection{Decay}
\subsubsection{Convergence to the attractor}
\begin{figure}[t]
 \begin{centering}
   \psfrag{t}{$\tau$} {\psfrag{phi}{$\|\Phi-\Phi_{(a,b)}\|_{L_2}$}
     \includegraphics[width=0.41\textwidth]{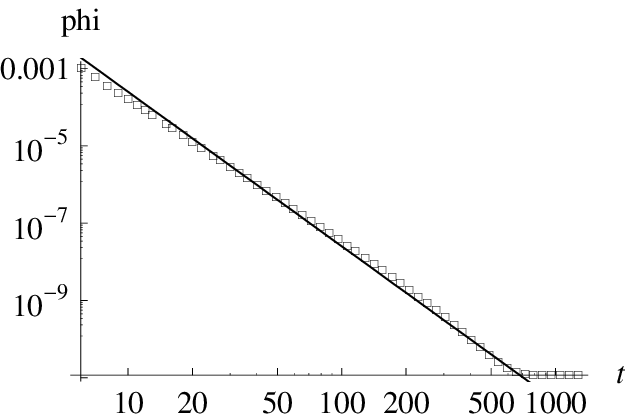}}\hspace{13mm}
   \psfrag{delta}{$\delta$}
   \includegraphics[width=0.41\textwidth]{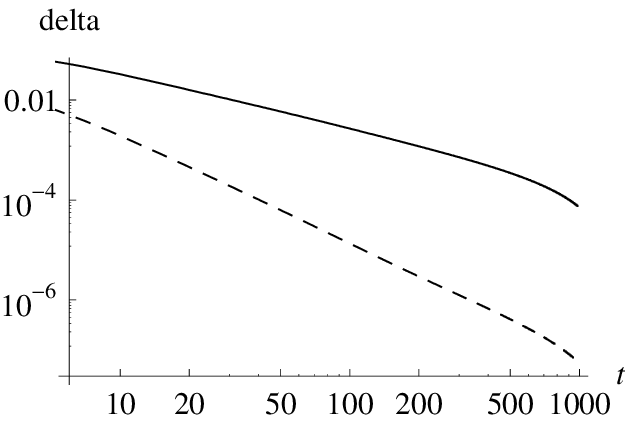}
   \caption{Left panel: The $L_2$-norm of the difference between a
     generic decaying numerical solution and the fitted attractor
     solution is depicted by the small squares on a log-log
     scale. The solid line, shown for comparison, has the slope $-4$, in accordance with Conjecture~1.
     Right panel: Modulation in the relative errors
     $\delta$ of the parameters $a$ (solid line) and $b$ (dashed
     line) on a log-log scale. The modulation for $a$ seems to decay
     as $t^{-1}$ and for $b$ as $t^{-2}$. The relative error, say
     $\delta_a$ for $a$, is defined as $\delta_a=|a(t)-a|/a$ where
     $a$ is the value of the parameter at the last time
     step of the numerical evolution. \label{fig:decay}}
 \end{centering}
\end{figure}

According to Conjecture~1, the difference between the attractor solution (with optimal
parameters) and a generic solution should decay as $t^{-4}$ for small data (\ref{conver}). This
behavior is confirmed numerically in Figure \ref{fig:decay}.

For this plot, we first fit the numerical solution to the attractor (\ref{phiab_hypal}) to
determine the optimal parameters $(a,b)$. We applied two methods for the fit. Fit in time, and
fit in space. For the fit in time we fix a grid point $\rho=\rho_f$, fit the solution
$\Phi(\tau,\rho_f)$ in $\tau$ and repeat the process for each grid point $\rho_f$. This gives us
a set of values $(a,b)_{\rho_f}$ depending on $\rho_f$. The average of these values over all grid
points gives the optimal parameters. For this method, the starting time for the fit needs to be
taken well behind the approximate wavefront because the convergence of the numerical solution to
the attractor family, and therefore the variation of $(a,b)_{\rho_f}$ over the grid, will depend
strongly on the foliation in early times. Note that our grid has infinite physical extent due to
compactification. The variation is expected to decrease in time which we have checked
numerically.  Choosing a large value of $\tau$ as the starting time of the time fit ensures
higher accuracy in the determination of the optimal parameters.

The second method that we applied is the fit in space in which we fix
a time $\tau=\tau_f$, fit the solution $\Phi(\tau_f,\rho)$ in $\rho$
and repeat the process on each time surface $\tau_f$. This gives us a
set of values $(a,b)_{\tau_f}$ depending on $\tau_f$. The resulting
time variation of the parameters is commonly referred to as
modulation. If our conjecture is correct, the modulation should be
small.

At the end, both methods should deliver the same values up to a small
error. This is used as a cross-check for the quality of fitting. The
modulation of the parameters $a$ and $b$ in time is shown on the right
panel in Figure \ref{fig:decay}. We observe that the accuracy of the
fitting in $b$ is better than in $a$. This is due to the fact that the
perturbation generated by changing $a$ decays in time while the
perturbation generated by changing $b$ does not (\ref{gena}), hence an
error in determining $b$ is easier to spot.

Once the optimal parameters have been determined, we compute the difference between the attractor
solution and the numerical solution at each time step. The difference $\Phi-\Phi_{(a,b)}$ at late
times has a
 strict sign for all values of $\rho$ on our numerical
 grid. Therefore the $L_2$-norm of this difference over the grid
 as a function of time gives a strong uniform measure of how close the
 solutions are to each other. Figure \ref{fig:decay} indicates that this
difference falls off as $t^{-4}$, as claimed in Conjecture~1.

\subsubsection{Exceptional solutions}
\begin{figure}[t]
 \begin{centering}
   \psfrag{t}{$\tau$} \psfrag{p}{$p_\rho$}
     \includegraphics[width=0.37\textwidth]{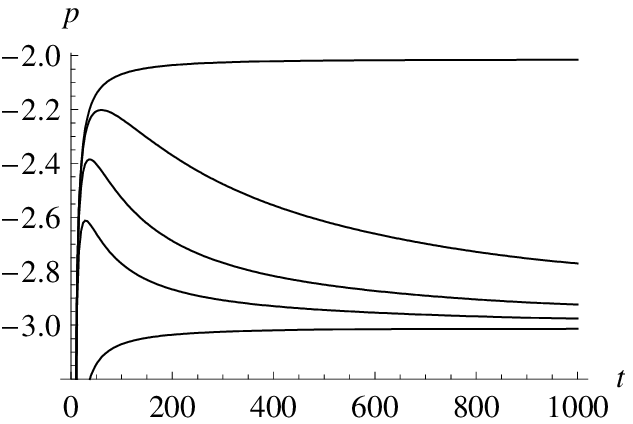}\hspace{13mm}
     \includegraphics[width=0.37\textwidth]{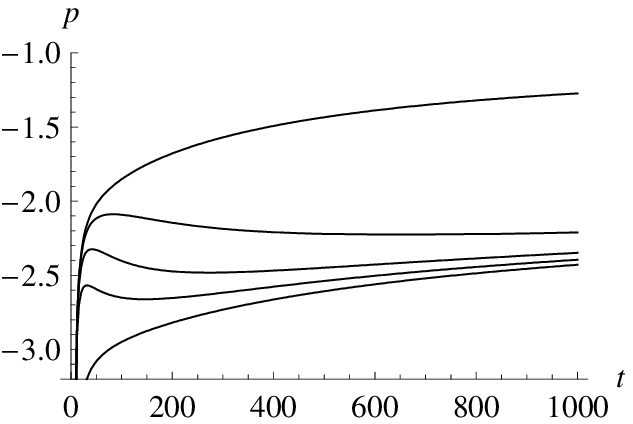}
     \caption{On the left panel we plot decay rates for the flip
       solution with $A_f=2.4913$ along the surfaces $\scri^+$,
       \mbox{$r=\{100,33,14,0\}$} from top to bottom. On the right
       panel we plot the decay rates for a solution with initial
       amplitude $A_f-0.01$ along the same surfaces. Here, the
       generic decay rates are obtained after a much longer time than
       that in Figure \ref{fig:tests}. \label{fig:flip}}
 \end{centering}
\end{figure}

As mentioned above in Remark~3, not all globally regular solutions converge to the attractor
$\phi_{(a,b)}$.  The existence of exceptional solutions with different asymptotic behavior can be
seen as follows. Consider a one-parameter family of initial data, such as (\ref{eq:id}). It turns
out that along such a family typically there occurs a flip of sign of the attractor, that is,
there is a (subcritical) amplitude $A_f$, such that for $A<A_f$ the solutions converge to, say,
$\phi_{(a,b)}$ and for $A>A_f$ they converge to $-\phi_{(a,b)}$. Performing a bisection we may
easily fine-tune the amplitude to $A_f$ and generate a special solution (below referred to as the
flip solution) that is not an element of the attractor family. One can expect that the flip
solution will have a faster decay rate than the generic global solutions. This expectation is
confirmed in Figure \ref{fig:flip} which shows that the flip solution decays as $t^{-3}$ at
timelike infinity and as $t^{-2}$ along null infinity. By modifying the initial amplitude
slightly away from the flip solution, we can see that the decay rates of the generic solutions
are attained at late times after a transient phase. Note that by construction the flip solutions
correspond to codimension-one initial data. It is likely that there exist initial data of higher
codimensions which lead to globally regular solutions with even faster decay rates, however the
numerical construction of such solutions would be very difficult.

\subsection{Blowup}
\subsubsection{Convergence to the attractor}
Merle and Zaag proved in \cite{mz1} that the ODE solution $\sqrt{2}/(T-t)$ determines the
universal rate of blowup for  equation (\ref{eqo}), however the problem of profile of blowup
remains open.  Numerical simulations in spherical symmetry \cite{bct} showed that for a solution
which blows up at the origin, its deviation from $\sqrt{2}/(T-t)$ near the tip of the past light
cone of the blowup point is very well approximated by the second eigenmode in the expansion
(\ref{bct}) (see Figure 4 in \cite{bct})
\begin{equation}\label{bct2}
   \phi(t,r)-\frac{\sqrt{2}}{T-t} \approx c_1 \left(1-\frac{r^2}{(T-t)^2}\right)\,,
\end{equation}
however an error on the right hand side was not quantified in \cite{bct}. According to
Conjecture~2
the approximation (\ref{bct2}) may be improved to
\begin{equation}\label{conj22}
    \phi(t,r)-\phi_{(a,b)}(t,r) \approx  \mathcal{O}((T-t)^2)\,,
\end{equation}
provided that the parameters $(a,b)$ are optimal.

\begin{figure}[t]
 \begin{centering}
   \psfrag{del}{$T-t$} \psfrag{dif}{$\phi-\phi_{(a,b)}$}
   \includegraphics[width=0.46\textwidth]{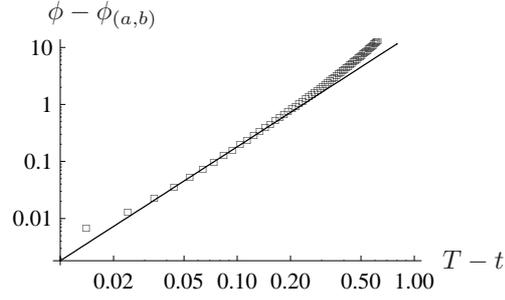}
     \caption{Blowup at the origin in standard coordinates on a
       log-log scale. The difference between the attractor solution
       and the numerical solution is depicted by small squares. The linear fit (solid line)
       gives the slope $2.02$, which confirms Conjecture~2.
       The deviation from the straight line seen at times very close
       to the blowup time is due to the fact that the domain of
       analysis extends beyond the past light cone of the blowup
       point.
\label{fig:st_blup}}
 \end{centering}
\end{figure}

The numerical verification of the formula (\ref{conj22}) is shown in Figure~\ref{fig:st_blup}.
The data for this plot were produced using a code based on standard Minkowski coordinates.

We point out that there is no genericity condition in Conjecture~2. This might appear surprising
in view of existence of a countable family of regular self-similar solutions of equation
(\ref{eq}) \cite{bbmw}. It seems that these self-similar solutions do not play any role in the
Cauchy evolution, which is probably due to the fact that they contain singularities outside the
past light cone of the blowup point (see section 6 in \cite{bbmw}).

\subsubsection{Blowup surface}
\begin{figure}[t]
 \begin{centering}
   \psfrag{r}{$\rho$} \psfrag{phi}{$\Phi$}
     \includegraphics[width=0.46\textwidth]{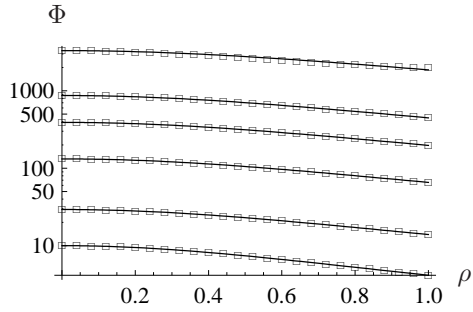}
     \caption{Simultaneous blowup along the numerical grid for fine
         tuned initial data. The solution is plotted on a log
       scale at various time steps close to the blowup time $T$. The
       solid lines depict an attractor solution with $b=-1/2$ at the
       corresponding times. The times are from bottom to top
       $\triangle T = T-\tau =
       \{0.34,0.1,0.02,0.007,0.003,0.0008\}$. \label{fig:globup}}
 \end{centering}
\end{figure}

For $b<0$ the attractor solutions blow up along a hyperboloid
\begin {equation}\label{blowatt}
t= -a -\frac{1}{2b} + \sqrt{\frac{1}{4b^2}+r^2}\,.
\end{equation}
This surface has the form (\ref{eq:tau}) with the mean extrinsic curvature $K=-6b$. Hence, for
our choice of the hyperboloidal foliation with $K=3$ the blowup surface of the attractor with
$b=-1/2$ coincides with one leaf of the foliation. Therefore, if Conjecture~2 is correct, in the
case of the blowup solution converging to the attractor with the optimal parameter $b=-1/2$ we
should observe an approximately simultaneous blowup along the whole grid. This expectation is
verified in Figure \ref{fig:globup}.

Note that, for a given constant mean curvature foliation, the simultaneous blowup is not generic.
To produce data for Figure \ref{fig:globup} we used the dependence of $b$ on the amplitude $A$ of
the gaussian. For relatively small (but supercritical) values of $A$ we have $-1/2<b<0$ and the
blowup first occurs at null infinity, while for larger amplitudes we have $b<-1/2$ and the blowup
first occurs at the origin. Performing bisection between these two states
we fine-tuned to the blowup solution with $b=-1/2$.

\begin{figure}[t]
 \begin{centering}
   \psfrag{rho}{$\rho$} \psfrag{phi}{$\Phi$}
     \includegraphics[width=0.46\textwidth]{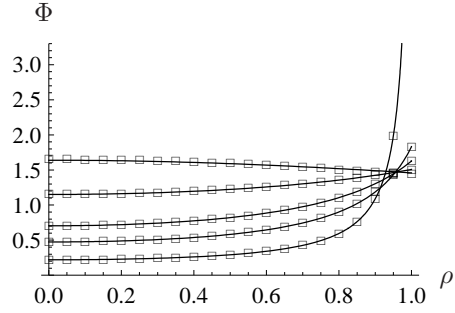}
     \caption{Blowup at null infinity. The numerical
       solution (small squares) at the times (counting, near the origin, from top to bottom)
       $\tau = \{2.3,3.2,5,7,25\}$ is compared to
       the attractor solution (solid line) with
       $b=-0.02$.  We see that the solution grows at
       null infinity while it decays near the
       origin. \label{fig:nullup}}
 \end{centering}
\end{figure}

We would like to emphasize that a rather counter-intuitive phenomenon of blowup at null infinity
is a mere coordinate effect. This only occurs if the mean extrinsic curvature of the blowup
 surface, given by $-6b$, is smaller then the mean extrinsic curvature
 of the hyperboloidal foliation, $K$, used in the numerical
 simulation.  Nevertheless, as discussed above, we can use this
effect to our advantage to probe the shape of the blowup surface. Note
also that one can predict the blowup at null infinity by fitting the
attractor to the numerical solution at the origin. If this fit gives a
value of $b \in (-1/2,0)$, then we know that the solution will blow up
at null infinity even though it is decaying at the origin (see Figure
\ref{fig:nullup}).

\subsection{Critical behavior}
Now we consider the behavior of solutions for initial data lying at the boundary between
dispersion and blowup. Let us recall that this problem was addressed before in \cite{bct} for the
focusing wave equation $\partial_{tt} \phi - \Delta\phi- \phi^p=0$ with three values of the
exponent $p$ (corresponding to three different criticality classes with respect to scaling of
energy): $p=3$ (subcritical), $p=5$ (critical), and $p=7$ (supercritical). It was shown there
that the nature of the critical solution, whose codimension-one stable manifold separates blowup
from dispersion, depends on $p$: for $p=7$ the critical solution is self-similar, while for $p=5$
it is static. For $p=3$ the critical solution was not determined because of numerical
difficulties (although with hindsight it could have been inferred from Figure~11 in \cite{bct}).
Now, in view of Conjectures~1 and 2 which assert that the solution $\phi_{(a,b)}$ is an attractor
for generic dispersive solutions if $b>0$ and for blowup solutions if $b<0$, it is easy to guess
that the critical solution corresponds to $b=0$, hence it has the form $\phi_c=\sqrt{2}/(t+a)$.

\begin{figure}[t]
 \begin{centering}
   {\psfrag{rho}{$\rho$}\psfrag{phi}{$\Phi$}
     \includegraphics[width=0.44\textwidth]{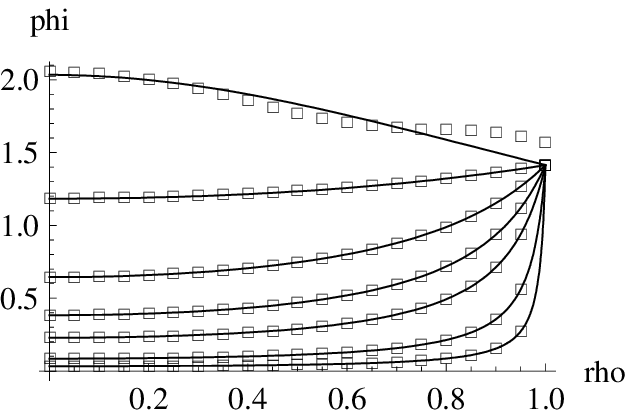}\hspace{13mm}}
   {\psfrag{t}{$\tau$} \psfrag{phi}{$\Phi|_{\scri^+}$}
     \includegraphics[width=0.44\textwidth]{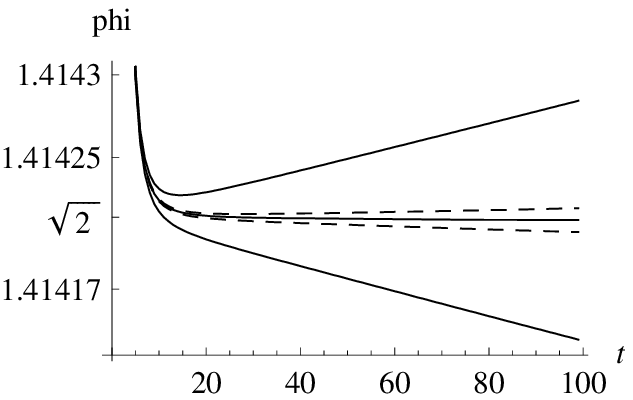}}
   \caption{Left panel: The critical solution on the grid at times
     $\tau=\{2,3,5,8,13,34,89\}$ (counting from top to bottom). We see
     that at $\tau=2$ the solution is not yet described well by the
     attractor solution, but already after $\tau=3$ the agreement is
     good. Right panel: Numerical solutions close to the critical
     solution evaluated at null infinity. Denoting the initial
     amplitude for the solution that corresponds to the almost constant solid line at $\sqrt{2}$ by
     $A_c$, the deviating dashed curves have $A=A_c \pm 10^{-9}$, the
     deviating solid ones have $A=A_c\pm 10^{-8}$. \label{fig:critical}}
 \end{centering}
\end{figure}

The critical solution is difficult to study with standard numerical methods by bisection because
it is a globally decaying solution, but it can be studied very accurately with the conformal
method. The reason is that the conformal scaling factors out the leading asymptotic
behavior implying that the $1/t$ decay of the critical solution is factored out at the conformal
boundary.  Specifically, the rescaled critical solution $\Phi_c=\phi_c/\Omega$ in the new
coordinates as given in (\ref{phiab_hypal}) evaluated at null infinity becomes
$\Phi_c|_{\scri^+}=\Phi_{(a,0)}(\tau,1) = \sqrt{2}$, hence the deviation of $\Phi|_{\scri^+}$
from $\sqrt{2}$ can be used as a bisection criterion. On the left panel of Figure
\ref{fig:critical} we plot the critical solution on the grid at various time steps and compare it
to the attractor solution with $b=0$ and fitted $a$. We see that the solution is decaying while
its value at $\scri^+=\{\rho=1\}$ is constant. The instability of the critical solution can be seen
by evolving initial data that differ slightly from the critical amplitude $A_c$ -- for such data
the deviation from $\sqrt{2}$ at $\scri^+$ grows linearly with time. This is depicted on the
right panel of Figure \ref{fig:critical} for four different values of marginally critical
amplitude.

\section{Final remarks}

We wish to emphasize that the use of the hyperboloidal foliation (\ref{eq:tau}) in combination with
the conformal method was instrumental in unraveling the dynamics of global solutions of
equation (\ref{eq}). First and foremost, this method eliminates the need of introducing an
artificial boundary, which is a notorious problem in computing wave propagation on unbounded
domains. Second, the intersection of $\tau=const$ hyperboloids with $\scri^+$ increases
monotonically with $\tau$ leading to the dispersive dissipation of energy along the leaves of the
foliation which is a mechanism responsible for convergence to the attractor. Third,
the conformal rescaling
allows a very accurate computation of the critical solution by factoring out its leading
asymptotic behavior. Finally, with this approach we can probe efficiently the shape of the blowup
surface and observe the simultaneous blowup on the whole grid.

Our main observation that a simple family of exact solutions can act as a universal attractor for
solutions of the nonlinear wave equation was unexpected to us. It is clear that this surprising
phenomenon is intimately related to the conformal invariance of the cubic wave equation, and
therefore it is more a curiosity than a stable property, in particular it is absent for
semilinear focusing wave equations
\begin{equation}\label{focusp}
\Box \phi +|\phi|^{p-1}\phi=0
\end{equation}
with $p\neq 3$. However, we conjecture that the threshold behavior
mediated by the slowly decaying global solution is structurally stable in the sense that the ODE
solution of equation~(\ref{focusp})
\begin{equation}\label{diffp}
   \phi=\frac{c}{(t+a)^{\frac{2}{p-1}}}\,,\quad \mbox{where} \qquad
   c=\left[\frac{2(p+1)}{(p-1)^2}\right]^{\frac{1}{p-1}}\,,
\end{equation}
is a critical solution for all exponents satisfying $1+\sqrt{2}<p\leq 3$. Note that the decay
rate exponent of the critical solution, $2/(1-p)$, and the decay rate exponent of generic
solutions, $1-p$, merge for $p\rightarrow 1+\sqrt{2}$, which is consistent with the fact that for
$p \leq 1+\sqrt{2}$ all solutions of equation (\ref{focusp}) with compactly supported initial
data blow up in finite time \cite{fjohn}.

\ack We thank Helmut Friedrich, Sascha Husa, Vince Moncrief, and Nikodem Szpak for discussions.
PB thanks the Albert Einstein Institute in Golm for hospitality and support during the initial
phase of this project. PB acknowledges support by the MNII grants: NN202 079235 and
189/6.PRUE/2007/7. AZ acknowledges support by the NSF grant PHY0801213 to the University of
Maryland.

\section*{References}
\bibliography{references}\bibliographystyle{utphys}

\providecommand{\href}[2]{#2}\begingroup\raggedright\begin{thebibliography}{10}

\bibitem{ch1}
D.~Christodoulou, ``Global solutions of nonlinear hyperbolic equations for
  small initial data,'' {\em Comm. Pure Appl. Math.} {\bf 39} (1986) no.~2,
  267--282.

\bibitem{nik}
N.~Szpak, P.~Bizo\'n, T.~Chmaj, and A.~Rostworowski, ``Linear and nonlinear
  tails ii: exact decay rates in spherical symmetry,'' {\em J. Hyperbolic Diff.
  Equations} {\bf 6} (2009)  107--125,
  \href{http://arxiv.org/abs/arXiv:0712.0493}{{\tt arXiv:0712.0493}}.

\bibitem{mz1}
F.~Merle and H.~Zaag, ``{Determination of the blow-up rate for a critical
  semilinear wave equation},'' {\em Math. Ann.} {\bf 331} (2005)  395--416.

\bibitem{bct}
P.~Bizo\'n, T.~Chmaj, and Z.~aw~Tabor, ``On blowup for semilinear wave
  equations with a focusing nonlinearity,'' {\em Nonlinearity} {\bf 17} (2004)
  2187--2201, \href{http://arxiv.org/abs/math-ph/0311019}{{\tt
  math-ph/0311019}}.

\bibitem{anco}
S.~C. Anco and S.~Liu, ``Exact solutions of semilinear radial wave equations in
  n dimensions,'' {\em Journal of Mathematical Analysis and Applications} {\bf
  297} (2004)  317.

\bibitem{ch2}
D.~Christodoulou, {\em Mathematical Problems in General Relativity I, (Zurich
  Lectures in Advanced Mathematics)}.
\newblock European Mathematical Society, Switzerland, 2008.

\bibitem{gp}
V.~Galaktionov and S.~Pohozaev, ``{On similarity solutions and blow-up spectra
  for a semilinear wave equation},'' {\em Quart. Appl. Math.} {\bf 61} (2003)
  583--600.

\bibitem{Penrose65}
R.~Penrose, ``Zero rest-mass fields including gravitation: Asymptotic
  behaviour,'' {\em Proc. Roy. Soc. Lond.} {\bf A284} (1965)  159--203.

\bibitem{Hawking73}
S.~W. Hawking and G.~F.~R. Ellis, {\em The large scale structure of spacetime}.
\newblock Cambridge University Press, Cambridge, England, 1973.

\bibitem{Zeng07}
A.~Zengino\u{g}lu, ``{Hyperboloidal foliations and scri-fixing},''
  \href{http://dx.doi.org/10.1088/0264-9381/25/14/145002}{{\em Class. Quant.
  Grav.} {\bf 25} (2008)  145002},
\href{http://arxiv.org/abs/0712.4333}{{\tt arXiv:0712.4333 [gr-qc]}}.

\bibitem{Moncrief00}
V.~Moncrief, ``Conformally regular {ADM} evolution equations,,'' 2000.
\newblock Talk at Santa Barbara, \\
  \texttt{http://online.itp.ucsb.edu/online/numrel00/moncrief}.

\bibitem{Husa02}
S.~Husa, ``{Numerical relativity with the conformal field equations},'' {\em
  Lect. Notes Phys.} {\bf 617} (2003)  159--192,
\href{http://arxiv.org/abs/gr-qc/0204057}{{\tt arXiv:gr-qc/0204057}}.

\bibitem{Fodor03}
G.~Fodor and I.~Racz, ``{What does a strongly excited 't Hooft-Polyakov
  magnetic monopole do?},''
  \href{http://dx.doi.org/10.1103/PhysRevLett.92.151801}{{\em Phys. Rev. Lett.}
  {\bf 92} (2004)  151801},
\href{http://arxiv.org/abs/hep-th/0311061}{{\tt arXiv:hep-th/0311061}}.

\bibitem{bbmw}
P.~Bizo\'n, P.~Breitenlohner, D.~Maison, and A.~Wasserman, ``Self-similar
  solutions of the cubic wave equation,''
  \href{http://arxiv.org/abs/arXiv:0712.0493}{{\tt arXiv:0712.0493}}.

\bibitem{fjohn}
F.~John, ``Blow-up of solutions of nonlinear wave equations in three space
  dimensions,'' {\em Manuscript. Math.} {\bf 28} (1979)  235--268.

\end{thebibliography}\endgroup

\end{document}